\documentclass[12pt]{elsarticle}

\usepackage{amssymb}
\usepackage{xcolor}
\usepackage{mathrsfs}
\usepackage{amscd}
\usepackage{appendix}
\usepackage{geometry}

\usepackage{amsfonts}
\usepackage{amsmath}
\newtheorem{theorem}{Theorem}[section] % 1st argument is your name for it
\newtheorem{lemma}[theorem]{Lemma}     % 2nd argument is what is printed

\newtheorem{remark}[theorem]{Remark}

\numberwithin{equation}{section}

\def\Proof{
\noindent \it Proof.\ \ \rm}
\def\qedbox{$\rlap{$\sqcap$}\sqcup$}
\journal{arXiv.com}

\begin{document}
\newgeometry{left=1.9cm,right=2.2cm,top=2.0cm,bottom=2.5cm}
\begin{frontmatter}



\title{Log-H\"older Continuity of the Lyapunov Exponent for Jacobi Operators with Potentials Given by the Skew-Shift}

\author[1]{Licheng Fang}
\ead{flc@stu.ouc.edu.cn}

\author[1]{Daxiong Piao \fnref{b}}
\ead{dxpiao@ouc.edu.cn}
\address[1]{School of Mathematical Sciences,  Ocean University of China,
Qingdao
266100, P.R.China}
\fntext[b]{Corresponding author}
\begin{abstract}
 In this paper we study one-dimensional Jacobi operators on the lattice with a potential given by the skew shift. We show that the large deviation theorem takes place for Diophantine frequency and sufficiently large disorder. Combining the large deviation theorem with the avalanche principle, we prove the log-H\"older continuity of the Lyapunov exponent.
\end{abstract}


\end{frontmatter}
\section{Introduction}
We consider the quasi-periodic Jacobi operators with skew-shift on $\ell^{2}(\mathbb{Z})$ defined by
$$(H_{\omega,(x,y)}\phi)(n)=-a(y+(n+1)\omega)\phi(n+1)-a(y+n\omega)\phi(n-1)+\lambda v(T^{n}_{\omega}(x,y))\phi(n)$$
where $T_{\omega}: \mathbb{T}^{2}\rightarrow \mathbb{T}^{2} $, $T_{\omega}(x,y)=(x+y,y+\omega)$, $a:\mathbb{T} \rightarrow \mathbb{R}$, $v: \mathbb{T}^{2} \rightarrow \mathbb{R}$ are real analytic function, and $1\leq|a(y)|\leq2$, $\omega\in(0,1)$ is the frequency and satisfies Diophantine condition. Specifically, we have
\begin{equation}\label{Diophantine}
\omega\in\Omega_{\varepsilon}:=\{\omega:\lVert n\omega\rVert>\frac{\varepsilon}{ n(\log n)^{2}}\;\; {\rm for \; any} \; n\in \mathbb{Z}^{+}\}
\end{equation}
with $\varepsilon\ll 1$. It is clear that
$\textrm{mes}[\mathbb{T}\setminus\Omega_{\varepsilon}]<C\varepsilon$
with an absolute constant $C$.

The Jacobi operator in fact is a second order, symmetric difference expression
\begin{align*}
H: \;\;\;\;\ell(\mathbb{Z}^{2})&\mapsto\ell(\mathbb{Z}^{2}) \\
\phi(n)&\mapsto-a(y+(n+1)\omega)\phi(n+1)-a(y+n\omega)\phi(n-1)+\lambda v(T^{n}_{\omega}(x,y))\phi(n)
 \end{align*}
It is associated with the tridiagonal matrix $(v_{n}=v(T_{\omega}^{n}(x,y)), a_{n}=a(y+n\omega))$
\[\left( {\begin{array}{*{20}{c}}
\ddots&\ddots&\ddots\\
-a_{n}&\lambda v_{n}&-a_{n+1}\\
&-a_{n+1}&\lambda v_{n+1}&-a_{n+2}\\
&&-a_{n+2}&\lambda v_{n+2}&-a_{n+3} \\
&&\ddots&\ddots&\ddots
\end{array}} \right)\]

Here we introduce two difference expressions
$$(\partial \psi)(n)=\psi(n+1)-\psi(n),\;\;\; (\partial^{*}\psi)(n)=\psi(n-1)-\psi(n),$$
using this we can rewrite $H$ in the following way
\begin{align}\label{eq2}
(H\psi)(n)=(\partial a\partial^{*}\psi)(n)+(\lambda v(T^{n}_{\omega}(x,y))-a_{n}-a_{n+1})\psi(n).
\end{align}
By means of the summation by parts formula (also known as Abel transform)
\begin{equation}\label{eq3}
\sum_{j=m}^{n}\psi(j)(\partial \varphi)(j)=\psi(n)\varphi(n+1)-\psi(m-1)\varphi(m)+\sum_{j=m}^{n}(\partial^{*}\psi)(j)\varphi(j),
\end{equation}
we have the Green's formula for discrete case
\begin{equation}\label{eq51}
\sum_{j=m}^{n}\big(\varphi(H \psi)-(H \varphi)\psi\big)(j)=W_{n}(\varphi,\psi)-W_{m-1}(\varphi,\psi),
\end{equation}
where $W_{n}$ is the modified Wronskian
$$W_{n}(\varphi,\psi)=a_{n+1}(\psi(n)\varphi(n+1)-\varphi(n)\psi(n+1)).$$

Now we considering the Jacobi difference equation
\begin{equation}\label{eq4}
H\phi=E\phi.
\end{equation}
If $\varphi$ and $\psi$ are both solution of  (\ref{eq4}), the left side of equation (\ref{eq51}) is zero, which means the Wronskian is constant (i.e., does not depend on $n$). We always omit the index $n$ for this case. Moreover, it is nonzero if and only if $\varphi$ and $\psi$ are linearly independent.
Since the space of solution is two dimensional, we can pick two linearly independent solution $\varphi$, $\psi$ and write any solution $\phi$ of (\ref{eq4}) as a linear combination of these two solutions
\begin{equation}\label{eq5}
\phi(n)=\frac{W(\phi,\psi)}{W(\varphi,\psi)}\varphi(n)-\frac{W(\phi,\varphi)}{W(\varphi,\psi)}\psi(n).
\end{equation}

For further studying, it is convenient to introduce the following fundamental solutions $\varphi$, $\psi$ $\in \ell^{2}(\mathbb{Z})$
$$H \varphi(E,\cdot,n_{0})=E\varphi(E,\cdot,n_{0}),\;\;\; H \psi(E,.,n_{0})=E\psi(E,\cdot,n_{0}),$$
fulling the initial conditions
$$\varphi(E,n_{0},n_{0})=1,\;\;\;\varphi(E,n_{0}+1,n_{0})=0,$$
$$\psi(E,n_{0},n_{0})=0,\;\;\;\psi(E,n_{0}+1,n_{0})=1.$$
Since the Wronskian of $\varphi(E,\cdot,n_{0})$ and $\psi(E,\cdot,n_{0})$ does not depend on $n$ we can evaluate it at $n_{0}$
$$W(\varphi(E,\cdot,n_{0}),\psi(E,\cdot,n_{0}))=-a_{n_{0}+1}$$
and consequently equation (\ref{eq5}) is simplified to
\begin{equation}\label{eq6}
\phi(n)=\phi(n_{0})\varphi(E,n,n_{0})+\phi(n_{0}+1)\psi(E,n,n_{0}).
\end{equation}

Eigenvalue equation
 $$ -a_{n+1}\phi(n+1)-a_{n}\phi(n-1)+\lambda v(T^{n}_{\omega}(x,y))\phi(n)=E\phi(n)$$
can be rewritten as
$$\Big(\begin{array}{ccc}\phi(n+1)\\ \phi(n) \end{array}\Big)
 =\frac{1}{a_{n+1}}\Big(\begin{array}{ccc} \lambda v(T_{\omega}^{n}(x,y))-E & -a_{n} \\ a_{n+1}  & 0 \end{array} \Big)\Big(\begin{array}{ccc}\phi(n)\\ \phi(n-1) \end{array}\Big)=A_{n}(x,y;\lambda,E)\Big(\begin{array}{ccc}\phi(n)\\ \phi(n-1) \end{array}\Big),$$
 or
 $$\Big(\begin{array}{ccc}\phi(n)\\ \phi(n-1) \end{array}\Big)
 =\frac{a_{n+1}}{a_{n}}\Big(\begin{array}{ccc}  0 & a_{n} \\ -a_{n+1}  & \lambda v(T_{\omega}^{n}(x,y))-E \end{array} \Big)\Big(\begin{array}{ccc}\phi(n+1)\\ \phi(n) \end{array}\Big)=A_{n}^{-1}(x,y;\lambda,E)\Big(\begin{array}{ccc}\phi(n+1)\\ \phi(n) \end{array}\Big).$$
 The matrix $A_{n}(x,y;\lambda,E)$ is often referred to as transfer matrix. The corresponding (non-autonomous) flow is given by the fundamental matrix
\begin{align}\label{eq7}
M_{[n,n_{0}]}(x,y;\lambda,E)&=\Big( \begin{array}{ccc}\psi(E,n+1,n_{0})&\varphi(E,n+1,n_{0})\\\psi(E,n,n_{0})&\varphi(E,n,n_{0})\end{array}\Big)\notag\\
&=\left\{
\begin{array}{rcl}
A_{n}(x,y;\lambda,E)\cdots A_{n_{0}+1}(x,y;\lambda,E)& &{n>n_{0}}\\
I\ \ \ \ \ \ \ \ \ \ \ \ \ \ \ \ \ \ \ \ \ \  & &{n=n_{0}}\\
 A_{n+1}^{-1}(x,y;\lambda,E)\cdots A_{n_{0}}^{-1}(x,y;\lambda,E)& & {n<n_{0}}
\end{array} \right.
\end{align}
More explicitly, equation (\ref{eq6}) is now equivalent to
$$\Big(\begin{array}{ccc}\phi(n+1)\\\phi(n)\end{array} \Big)=M_{[n,n_{0}]}(x,y;\lambda,E)\Big(\begin{array}{ccc}\phi(n_{0}+1)\\\phi(n_{0})\end{array} \Big)$$
and $M_{[n,n_{0}]}(x,y;\lambda,E)$ satisfies the usual group law
$$M_{[n,n_{0}]}(x,y;\lambda,E)=M_{[n,n_{1}]}(x,y;\lambda,E)M_{[n_{1},n_{0}]}(x,y;\lambda,E)$$
and constancy of the Wronskian implies
$$\det M_{[n,n_{0}]}(x,y;\lambda,E)=\frac{a_{n_{0}+1}}{a_{n+1}}.$$

Let's use $M_{[n,0]}(x,y;\lambda,E)=M_{n}(x,y;\lambda,E)$ and define the Lyapunov exponent
$$L(E):=\lim_{n\rightarrow\infty}L_{n}(E)=\inf_{n\geq 1}L_{n}(E),$$
where
$$L_{n}(E)=\int_{\mathbb{T}^{2}}\frac{1}{n}\log\lVert M_{[n,n_{0}]}(x,y;\lambda,E)\rVert dxdy.$$
This existence of this limit is guaranteed by subadditivity. By virtue of
$$\lVert M_{n_{0}}(x,y;\lambda,E)\rVert^{-1}\lVert M_{n}(x,y;\lambda,E)\rVert\leq\lVert M_{[n,n_{0}]}(x,y;\lambda,E)\rVert\leq\lVert M_{n_{0}}^{-1}(x,y;\lambda,E)\rVert\lVert M_{n}(x,y;\lambda,E)\rVert$$
we can find that the definition of $L(E)$ is indeed independent of $n_{0}$.

Next, using the property $M_{[n_{0},n_{1}]}(x,y;\lambda,E)=M_{[n_{1},n_{0}]}^{-1}(x,y;\lambda,E)$, we get the following
$$\Big(\begin{array}{ccc} \psi(E,n_{0}+1,n_{1})& \varphi(E,n_{0}+1,n_{1})\\\psi(E,n_{0},n_{1})&\varphi(E,n_{0},n_{1})\end{array}\Big)=\frac{a_{n_{1}+1}}{a_{n_{0}+1}}\Big(\begin{array}{ccc} \varphi(E,n_{1},n_{0})& -\varphi(E,n_{1}+1,n_{0})\\-\psi(E,n_{1},n_{0})&\psi(E,n_{1}+1,n_{0})\end{array}\Big)$$
by using (\ref{eq6}), and a straight calculation yields
\begin{equation*}
\psi(E,n,n_{0}+1)=-\frac{a_{n_{0}+2}}{a_{n_{0}+1}}\varphi(E,n,n_{0}),\ \ \psi(E,n,n_{0}-1)=\varphi(E,n,n_{0})+\frac{\lambda v_{n_{0}}-E}{a_{n_{0}+1}}\psi(E,n,n_{0}).
\end{equation*}
Let $J_{n_{0},n}$ be the Jacobi matrix
\[J_{n_{0},n} =\left( {\begin{array}{*{20}{c}}
\lambda v_{n_{0}+1}&-a_{n_{0}+2}\\
-a_{n_{0}+2}&\lambda v_{n_{0}+2}&\ddots\\
&\ddots&\ddots&\ddots\\
&&\ddots&\lambda v_{n-2}&-a_{n-1} \\
&&&-a_{n-1}&\lambda v_{n-1}
\end{array}} \right)\]
We notice that if $E$ is a zero of $\psi(\cdot,n,n_{0})$, then $(\psi(E,n_{0}+1,n_{0}),\cdots,\psi(E,n-1,n_{0}))$ is an eigenvector of $J_{n_{0},n}$ corresponding to the eigenvalue $E$. Since the converse statement is true, the polynomials (in $E$) $\psi(E,n,n_{0})$ and $\det (J_{n_{0},n}-E I)$ only differ by a constant which can be deduced from (\ref{eq7}). Hence we have the following expansion for $\psi(E,n,n_{0})$, $n>n_{0}$,
$$\psi(E,n,n_{0})=\frac{\det(J_{n_{0},n}-EI)}{\Pi_{j=n_{0}+2}^{n}a_{j}}.$$

So, if we define
\[f_{n}(x,y;\lambda,E) :=\det \left( {\begin{array}{*{20}{c}}
\lambda v_{1}-E&-a_{2}&0& 0&\cdots&\cdots &0\\
-a_{2}&\lambda v_{2}-E &-a_{3}&0&\cdots&\cdots&0\\
0&-a_{3}& \lambda v_{3}-E&-a_{4}&\cdots&\cdots&0\\
\vdots & \vdots &\vdots &\vdots &{}& \ddots &-a_{n}\\
0&0&{0}&{\cdots}&0&-a_{n}&\lambda v_{n}-E
\end{array}} \right)\]
then the fundamental matrix for Jacobi operator is
\begin{equation}
 M_{n}(x,y;\lambda,E)=\Bigg(\begin{array}{ccc}   \frac{f_{n}(x,y;\lambda,E)}{\prod_{j=2}^{n+1}a_{j}} & -\frac{a_{1}}
{a_{2}}\frac{f_{n-1}(T(x,y);\lambda,E)}{\prod_{j=3}^{n+1}a_{j}}\\
\frac{ f_{n-1}(x,y;\lambda,E)}{\prod_{j=2}^{n}a_{j}}&-\frac{a_{1}}{a_{2}}\frac{f_{n-2}(T(x,y);\lambda,E)}
{\prod_{j=3}^{n}a_{j}}\end{array} \Bigg)
\end{equation}

\subsection{Background and the main results}

When $a(y)\equiv1$, it is called Schr\"odinger operators, there are lots of continuity results related to this. Furthermore when the $T_{\omega}$ is the shift on the one-dimensional torus $\mathbb{T}$. You \cite{YZ} proved that $L(E)$ is H\"older continuity in $E\in[E_{1},E_{2}]$ when $L(E)>\gamma>0$ for $E\in[E_{1},E_{2}]$ and $\omega$ is weaker Liouville. Wang \cite{WY} proved that the discontinuity of $L(E)$ when the potential function is nonanalytic. In the general irrational frequency case, Bourgain \cite{BJ} proved that the continuity of $L(E)$ but without giving the specific regularity property. When $T_{\omega}$ is the skew-shift on the two-dimensional torus $\mathbb{T}^{2}$. Bourgain \cite{BGS1} proved that $L(E)$  is log-H\"older continuity for Diophantine frequency.

Apparently, the Jacobi operator is a more complex operator. Furthermore when the $T_{\omega}$ is the shift on the one-dimensional torus $\mathbb{T}$. Tao \cite{T}\cite{T1} proved that the $L(E)$ is H\"older continuity when frequency $\omega$ is weaker Liouville. In this paper we consider Jacobi operators case with potentials given by skew-shift, and we get the  log-H\"older continuity of $L(E)$ for Diophantine frequency. So in this case, this is a partial improved results of \cite{BGS1}.

In order to get the continuity of the Lyapunov exponent, we need to get the large deviation theorem that is the main section in this paper. And combing it with the avalanche principle (stated later) to get the desired result. However the avalanche principle cannot be applied to $M_{n}$, because it's possible that $|\det M_{n}|\nleqslant 1$. To work around this issue it is natural to use the following two matrixes associated with $M_{n}$:
$$M_{n}^{u}(x,y;\lambda,E):=\frac{M_{n}(x,y;\lambda,E)}{|\det M_{n}(x,y;\lambda,E)|^{\frac{1}{2}}},$$
$$M_{n}^{a}(x,y;\lambda,E):=A_{n}^{\prime}(x,y;\lambda,E)\cdots A_{1}^{\prime}(x,y;\lambda,E),$$
where
$$A_{n}^{\prime}(x,y;\lambda,E)=\Big(\begin{array}{ccc} \lambda v(T_{\omega}^{n}(x,y))-E & -a_{n} \\ a_{n+1}  & 0 \end{array} \Big).$$
Based on the definition, it is straightforward to check that
\begin{equation}
\log \lVert M_{n}^{a}(x,y;\lambda,E)\rVert=\log \lVert M_{n}(x,y;\lambda,E)\rVert+\sum\limits_{j=1}^{n}{\log|a_{j+1}|},
\end{equation}
\begin{equation}
\log \lVert M_{n}^{u}(x,y;\lambda,E)\rVert=\log \lVert M_{n}(x,y;\lambda,E)\rVert-\frac{1}{2}\log|\frac{a_{1}}{a_{n+1}}|,
\end{equation}
\begin{equation}\label{eq55}M_{n}^{a}(x,y;\lambda,E)=\prod_{j=1}^{n}(a_{j}a_{j+1})^{\frac12}M_{n}^{u}(x,y;\lambda,E),\end{equation}
and then we define
\begin{equation}
L_{n}^{a}:=\int_{\mathbb{T}^{2}} \frac{1}{n}\log \lVert M_{n}^{a}(x,y;\lambda,E)\rVert dxdy,\ L_{n}^{u}:=\int_{\mathbb{T}^{2}} \frac{1}{n}\log \lVert M_{n}^{u}(x,y;\lambda,E)\rVert dxdy.
\end{equation}
\begin{remark}\label{remark1}
(1) Note that
$$\lVert A_{n}(x,y;\lambda,E)\rVert\leq \frac{C(\lambda,v,E,a)}{|a(y+(n+1)\omega)|},$$
Therefore
$$\frac{1}{n}\log \lVert M_{n}(x,y;\lambda,E)\rVert\leq \log C(\lambda,v,E,a)-\frac{1}{n}\sum_{i=1}^{n}\log|a(y+(i+1)\omega)|$$
We always suppose that $|E|\leq E_{0}$, where $E_{0}$ depends on $\lambda,v,a$. For that matter we suppress $E$ from the notations of some of the constants involved.\\
(2) $\log \lVert M_{n}^{u}(x,y;\lambda,E)\rVert\geq0$, since $M_{n}^{u}(x,y;\lambda,E)$ is unimodular.\\
(3)\begin{align*}
0&\leq \frac{1}{n}\log \lVert M_{n}^{u}(x,y;\lambda,E)\rVert=\frac{1}{n}\log \lVert M_{n}(x,y;\lambda,E)\rVert-\frac{1}{2n}\log|\frac{a(y+\omega)}{a(y+(n+1)\omega)}|\\
&\leq \log C(\lambda,v,a)-\frac{1}{n}\sum_{i=1}^{n}\log|a(y+(i+1)\omega)|-\frac{1}{2n}\log|\frac{a(y+\omega)}{a(y+(n+1)\omega)}|
\end{align*}
(4)It is well-known fact that if $a$ is analytic function which is not identically zero then $(\log|a|)$ is integrable. Set
$$D=\int_{\mathbb{T}}\log|a|dy$$
Threfore $$\int_{\mathbb{T}^{2}}\Big|\frac{1}{n}\log\lVert M_{n}^{u}(x,y;\lambda,E)\rVert \Big|dxdy=\int_{\mathbb{T}^{2}}\frac 1n\log \lVert M_{n}^{u}(x,y;\lambda,E)\rVert dxdy\leq C^{\prime}(\lambda,v,a)-D:=C^{\prime\prime}(\lambda,v,a)$$
Similarly$$\int_{\mathbb{T}^{2}}(\frac{1}{n}\log\lVert M_{n}^{u}(x,y;\lambda,E)\rVert )^{2}dxdy\leq \tilde{C}(\lambda,v,a)$$
\end{remark}

It is well-known that for the case of the shift, i.e., $T(x,y)=(x+\omega_{1},y+\omega_{2})$, one takes advantage of the fact $\frac 1n\log\lVert M_{n}(\cdot)\rVert$ is subharmonic on a neighborhood of $\mathbb{T}^{2}$ and that it is bounded uniformly in $n$ in that neighborhood, so there are some useful results: if we let
\begin{equation}\label{eq13}
u_{n}(z_{1},z_{2})=\frac{1}{n}\log\lVert M_{n}(z_{1},z_{2};\lambda,E)\rVert,
\end{equation}
then we can get the decay of the Fourier coefficients
\begin{equation}\label{eq8}
|\widehat{u_{n}}(k_{1},k_{2})|\leq \frac{C(\lambda)}{|k_{1}|+|k_{2}|+1},
\end{equation}
and almost invariance
\begin{equation}\label{eq9}
\sup_{(x,y)\in\mathbb{T}^{2}}\Big|\frac 1K \sum_{k=1}^{K}u_{n}(T^{k}_{\omega}(x,y))-u_{n}(x,y)\Big|\leq C(\lambda)\frac{K}{n}.
\end{equation}

However, for the case of skew shift, i.e., $T(x,y)=(x+y,y+\omega)$, the subharmonic property is also needed and in this case, the width of holomorphicity in the second variable will have to be smaller than in the first by a factor of $\approx\frac{1}{n}$. This is due to the fact that the iteration of the skew shift is given by
$$T_{\omega}^{k}(x,y)=(x+ky+k(k-1)\omega/2,y+k\omega).$$
Complexifying in the variable $y$ therefore produces an imaginary part of size about n in half of the factors of the product $M_{n}$.
Considering that, here we introduce a scaling factor
\begin{equation}
S(\lambda,E)=\log(C_{v,a}+|\lambda|+|E|)\geq 1,
\end{equation}
 where $C_{v,a}$ is a constant depending on $v$, $a$ so that for all $n$,
\begin{equation}\label{eq10}
\sup\limits_{z_{1}\in\mathcal{A}_{\rho}}\sup\limits_{z_{2}\in\mathcal{A}_{\rho/n}}\frac{1}{n}\log\lVert M_{n}^{a}(z_{1},z_{2};\lambda,E)\rVert\leq S(\lambda,E).
\end{equation}

Here the main theorems are given
 \begin{theorem}\label{theorem1}
Assume $v$ is a nonconstant real analytic function on $\mathbb{T}^{2}$, let $\omega\in\Omega_{\varepsilon}$, then for all $\sigma<\frac{1}{24}$ there exist $\tau=\tau(\sigma)>0$ and constants $\lambda_{1}$ and $n_{0}$ depending only on $\varepsilon$, $v$ and $\sigma$ such that
\begin{equation}\label{eq31}
 \sup_{E} \textrm{mes}\Big[(x,y)\in\mathbb{T}^{2}:\big|\frac{1}{n}\log\lVert M_{n}(x,y;\lambda,E)\rVert-L_{n}(\lambda,E)\big|> S(\lambda,E)n^{-\tau}\Big]<C\exp(-n^{\sigma}),
 \end{equation}
for all $\lambda\geq \lambda_{1}$ and $n\geq n_{0}$. Furthermore, for those $\omega$, $v$ and all $E$,
$$L(\lambda,E)=\inf\limits_{n}L_{n}(\lambda,E)\geq \frac{1}{4}\log \lambda.$$
\end{theorem}
\begin{theorem}\label{theorem2}
Let $\omega$, $v$ and $\lambda_{1}$ be as in  Theorem \ref{theorem1}. For $\lambda>\lambda_{1}$, $L(E,\lambda)$ is log-H\"older continuous in $E$, that is
\begin{equation}
\Big|L(E,\lambda)-L(E^{\prime},\lambda)\Big| \leq C\exp(-c(\log|E-E^{\prime}|^{-1})^{\sigma}).
\end{equation}
\end{theorem}
\begin{remark}
There exists constant $E_{0}$ depending on $v, a$ such that the spectrum of the operator $H$ is contained in the interval $[-E_{0},E_{0}]$, and in the exterior of this interval, this system is uniformly hyperbolic, so we always suppose that $|E|\leq E_{0}$.
\end{remark}
\section{The main lemmas}
It is convenient that for a function $u$ defined on annulus $\mathcal{A}_{\rho}=\{z\in \mathbb{C}:1-\rho<|z|<1+\rho\}$, we shall write $u(z)$ instead of $u(re(x))\ (with \ e(x)=e^{2\pi ix})$. And for any positive integer $d$, $\mathbb{T}^{d}:=\mathbb{R}^{d}/\mathbb{Z}^{d}$ denotes the d-mimensional torus.
\begin{lemma}(\cite{GS} The Avalanche Principle )\label{Avalanche}
Let $A_{1},...,A_{n}$ be a sequence of $2\times2$-matrices. Suppose that
$$\max_{1\leq j\leq n}|\det A_{j}|\leq 1,$$
\begin{equation}\label{eq14}
\min_{1\leq j\leq n}\lVert A_{j}\rVert\geq\mu\geq n
\end{equation}
and
\begin{equation}\label{eq15}
\max_{1\leq j<n}[\log\lVert A_{j+1}\rVert+\log\lVert A_{j}\rVert-\log\lVert A_{j+1}A_{j}\rVert]\leq\frac12 \log\mu.
\end{equation}
Then
\begin{equation}\label{eq16}
\Big|\log\lVert A_{n},..., A_{1}\rVert+\sum_{j=2}^{n-1}\log\lVert A_{j}\rVert-\sum_{j=1}^{n-1}\log\lVert A_{j+1}A_{j}\rVert\Big|<C\frac{n}{\mu}.
\end{equation}
\end{lemma}
\begin{lemma}(\cite{BGS1})
Let $u:\mathbb{T}^{2}\rightarrow \mathbb{R}$ satisfy $\lVert u\rVert_{L^{\infty}(\mathbb{T}^{2})}\leq 1$. Assume that $u$ extends as a separately subharmonic function in each variable to a neighborhood of $\mathbb{T}^{2}$ such that for some $N\geq 1$ and $\rho>0$,
 $$\sup\limits_{z_{1}\in\mathcal{A}_{\rho}}\sup\limits_{z_{2}\in\mathcal{A}_{\rho}}|u(z_{1},z_{2})|\leq N.$$
 Furthermore, suppose that $u=u_{0}+u_{1}$ on $\mathbb{T}^{2}$ where
 \begin{equation}
 \lVert u_{0}-\langle u\rangle\rVert_{L^{\infty}(\mathbb{T}^{2})}\leq \varepsilon_{0}\;\;and\;\;\lVert u_{1}\rVert_{L^{1}(\mathbb{T}^{2})}\leq\varepsilon_{1}
 \end{equation}
 with $0<\varepsilon_{0},\;\varepsilon_{1}<1$. Here $\langle u\rangle:=\int_{\mathbb{T}^{2}}u(x,y)dxdy$. Then for any $\delta>0$,
\begin{equation} \textrm{mes}\Big[(x,y)\in\mathbb{T}^{2}: \big|u(x,y)-\langle u\rangle\big|>B^{\delta}  \log\frac{N}{\varepsilon_{1}}\Big]\leq CN^{2}\varepsilon_{1}^{-1}\exp(-cB^{-\frac 12+\delta}),\end{equation}
where $B=\varepsilon_{0}\log\frac{N}{\varepsilon_{1}}+N^{\frac{3}{2}}\varepsilon_{1}^{\frac 14}$. The constants $c$, $C$ only depend on $\rho$.
\end{lemma}
\begin{lemma}\label{lemma}(\cite{BGS1}) Let $u$ be 1-periodic subharmonic function defined on a neighborhood of $\mathbb{T}^{2}$. Suppose furthermore for some $\rho>0$,
$$\sup\limits_{z_{1}\in\mathcal{A}_{\rho}}\sup\limits_{z_{2}\in\mathcal{A}_{\rho}}|u(z_{1},z_{2})|\leq1.$$
For $\omega\in$ DC, any $\delta>0$ there exist constants $c,C$ depending on $\rho,\delta,\omega$ such that
\begin{equation}
\textrm{mes}\Big[(x,y)\in \mathbb{T}^{2}:\Big|\frac{1}{K}\sum_{k=1}^{K} u\circ T^{k}_{\omega}(x,y)-\langle u\rangle\Big|>K^{-\frac{1}{12}+2\delta}\Big]\leq C \exp(-cK^{\delta}),\end{equation}
for any positive integer $K$.
\end{lemma}

The following lemma provides the inductive step in the proof of the large deviation theorem.
\begin{lemma} \label{lemma1}
For $\omega\in\Omega_{\varepsilon}$, and suppose that $n,N$ are positive integers such that
\begin{equation}\label{eq11}
\textrm{mes}\Big[(x,y)\in\mathbb{T}^{2}:\Big|\frac{1}{n}\log\lVert M_{n}^{u}(x,y;\lambda,E)\rVert-L_{n}^{u}(\lambda,E)\Big|> S(\lambda,E)\frac{\gamma}{10}\Big]\leq N^{-10},
\end{equation}
\begin{equation}\label{eq12}
\textrm{mes}\Big[(x,y)\in\mathbb{T}^{2}:\Big|\frac{1}{2n}\log\lVert M_{2n}^{u}(x,y;\lambda,E)\rVert-L_{2n}^{u}(\lambda,E)\Big|> S(\lambda,E)\frac{\gamma}{10}\Big]\leq N^{-10} .
\end{equation}
And assume that
\begin{equation}
\min(L_{n}^{u}(\lambda,E),L_{2n}^{u}(\lambda,E))\geq\gamma S(\lambda,E),
\end{equation}
\begin{equation}
L_{n}^{u}(\lambda,E)-L_{2n}^{u}(\lambda,E)\leq\frac{\gamma}{40} S(\lambda,E),
\end{equation}
\begin{equation}\label{eq24}
9\gamma nS \geq10\log(2N) \;and \;n^{2} \leq N.
\end{equation}
Then there is some absolute constant $C_{0}$ with the property that
\begin{equation}
L_{N}^{u}(\lambda,E)\geq \gamma S(\lambda,E)-2(L_{n}^{u}(\lambda,E)-L_{2n}^{u}(\lambda,E))-C_{0}S(\lambda,E)nN^{-1}
\end{equation}
and
\begin{equation}\label{eq17}
L_{N}^{u}(\lambda,E)-L_{2N}^{u}(\lambda,E)\leq C_{0}S(\lambda,E)nN^{-1}.
\end{equation}
Further, for any $\sigma<\frac{1}{24}$ there is $\tau=\tau(\sigma)>0$ so that
\begin{equation}\label{eq21}
\textrm{mes}\Big[(x,y)\in\mathbb{T}^{2}:\Big|\frac{1}{N}\log\lVert M_{N}^{u}(x,y;\lambda,E)\rVert-L_{N}^{u}(\lambda,E)\Big|> S(\lambda,E)N^{-\tau}\Big]\leq C\exp(-N^{\sigma})
\end{equation}
with some constant $C=C(\sigma,\varepsilon)$.
\end{lemma}
\Proof
For simplicity, we write  $L_{N}^{u}=L_{N}^{u}(\lambda,E)$ etc. And we shall fix $\omega$, $\lambda$, and $E$. In particular, $S=S(\lambda,E)$. We write $N=mn+r,\ 0\leq r<n $. Denote the set on the left hand of (\ref{eq11}) by $\mathcal{B}_{n}$ and the set on the left-hand side of (\ref{eq12}) by $\mathcal{B}_{2n}$. For any $(x,y)\in\mathbb{T}^{2}/\mathcal{B}_{n}$,
\begin{equation}\label{eq40}\lVert M_{n}^{u}(x,y)\rVert\geq \exp(nL_{n}^{u}-nS\frac{\gamma}{10})\geq \exp(n\gamma S-nS\frac{\gamma}{10})=\exp(\frac{9\gamma}{10}nS)=:\mu>m.\end{equation}
For $1\leq j\leq m$, we consider $A_{j}=A_{j}(x,y):=M_{n}^{u}\circ T_{\omega}^{(j-1)n}(x,y)$. Then (\ref{eq40}) implies
\begin{equation}
\min_{1\leq j\leq m}\lVert A_{j}(x,y)\rVert\geq\mu\ \  {\rm for \ all}\  (x,y)\notin\bigcup\limits_{j=1}^{m}T_{\omega}^{-(j-1)n}\mathcal{B}_{n}.
\end{equation}
Furthermore, for any $(x,y)\notin\mathcal{B}_{0}:=\bigcup\limits_{j=1}^{m}T_{\omega}^{-(j-1)n}\mathcal{B}_{n}\cup\bigcup\limits_{j=1}^{m-1}T_{\omega}^{-(j-1)n}\mathcal{B}_{2n}$, with measure $<2N^{-10}\cdot N=2N^{-9}$, $1\leq j<m$, we have
\begin{align}
&\log\lVert A_{j+1}(x,y)\rVert+\log\lVert A_{j}(x,y)\rVert-\log\lVert A_{j+1}(x,y)\cdot A_{j}(x,y)\rVert\notag\\
&=\log\lVert M_{n}^{u}\circ T_{\omega}^{jn}(x,y)\rVert+\log\lVert M_{n}^{u}\circ T_{\omega}^{(j-1)n}(x,y)\rVert-\log\lVert M_{2n}^{u}\circ T_{\omega}^{(j-1)n}(x,y)\rVert\\
&\leq n(L_{n}^{u}+\frac{S}{10}\gamma)+n(L_{n}^{u}+\frac{S}{10}\gamma)-2n(L_{2n}^{u}-\frac{S}{10}\gamma)\notag\\
&\leq 2n(L_{n}^{u}-L_{2n}^{u})+\frac{4\gamma}{10}Sn\notag\\
&\leq \frac{9\gamma}{20}Sn=\frac 12\log \mu.
\end{align}
We can now apply the avalanche principle and get:
\begin{equation}\label{eq42}
\Big|\log\lVert A_{m}(x,y)\cdots A_{1}(x,y)\rVert +\sum_{j=2}^{m-1} \log\lVert A_{j}(x,y)\rVert -\sum_{j=1}^{m-1}\log\lVert A_{j+1}(x,y) A_{j}(x,y)\rVert \Big|\leq C\frac{m}{\mu} \end{equation}
for $(x,y)$ outside a set of measure $<2N^{-9}$. In particular, (\ref{eq42}) implies that
\begin{equation}
\Big|\log\lVert M_{mn}^{u}(x,y)\rVert +\sum_{j=2}^{m-1} \log\lVert M_{n}^{u}\circ T_{\omega}^{(j-1)n}(x,y) \rVert -\sum_{j=1}^{m-1}\log\lVert M_{2n}^{u}\circ T_{\omega}^{(j-1)n}(x,y)\rVert \Big|\leq C\frac{m}{\mu}.
\end{equation}
In view of  $N=mn+r$, then we have
\begin{equation}
\big|\log\lVert M_{N}^{u}(x,y)\rVert -\log\lVert M_{mn}^{u}(x,y)\rVert \big|\leq Sn\ \ and\  \ \big|\log\lVert M_{n}^{u}(x,y)\rVert\big|\leq Sn.
\end{equation}
Therefore
\begin{equation}\label{eq43}
\Big|\log\lVert M_{N}^{u}(x,y)\rVert +\sum_{j=0}^{m-1} \log\lVert M_{n}^{u}\circ T_{\omega}^{jn}(x,y) \rVert -\sum_{j=0}^{m-1}\log\lVert M_{2n}\circ T_{\omega}^{jn}(x,y)\rVert \Big|\leq C\frac{m}{\mu}+CSn\leq CSn.
\end{equation}
In (\ref{eq43}) replace $(x,y)$ by each of the elements $\{(x,y),T_{\omega}(x,y),\cdots,T_{\omega}^{n-1}\circ(x,y)\}$ and then average to get the following
\begin{equation}\label{eq44}
\Big|\frac{1} {N}\sum_{j=0}^{n-1}\frac{1}{n}\log\lVert M_{N}^{u}\circ T_{\omega}^{j}(x,y)\rVert +\frac{1}{N}\sum_{j=0}^{N-1} \frac{1}{n}\log\lVert M_{n}^{u}(x,y)\circ T_{\omega}^{j}(x,y) \rVert -\frac{2}{N}\sum_{j=0}^{N-1}\frac{1}{2n}\log\lVert M_{2n}^{u}\circ T_{\omega}^{j}(x,y)\rVert \Big|\leq C\frac{Sn}{N}.
\end{equation}
Due to the almost invariance property for skew-shift case
\begin{equation}\label{eq45} \big|\frac{1}{N}\log\lVert M_{N}^{u}(x,y)\rVert- \frac{1} {n}\sum_{j=0}^{n-1}\frac{1}{N}\log\lVert M_{N}^{u}\circ T_{\omega}^{j}(x,y)\rVert \big|\leq C\frac{Sn}{N}.\end{equation}
From (\ref{eq44}) (\ref{eq45}) we get
\begin{equation}
\Big|\frac{1}{N}\log\lVert M_{N}^{u}(x,y)\rVert +\frac{1}{N}\sum_{j=0}^{N-1} \frac{1}{n}\log\lVert M_{n}^{u}(x,y)\circ T_{\omega}^{j}(x,y) \rVert -\frac{2}{N}\sum_{j=0}^{N-1}\frac{1}{2n}\log\lVert M_{2n}^{u}\circ T_{\omega}^{j}(x,y)\rVert  \Big|\leq C\frac{Sn}{N}
\end{equation}
for $(x,y)\notin\mathcal{B}_{1}$, with measure $\textrm{mes}(\mathcal{B}_{1})\leq N\cdot\textrm{mes}(\mathcal{B}_{0})\leq 2N^{-8}$.
Integrating over $\mathbb{T}^{2}$ yields
\begin{equation}\big|L_{N}^{u}+L_{n}^{u}-2L_{2n}^{u} \big|\leq CSnN^{-1}+8SN^{-8}\leq C_{0}SnN^{-1}.\end{equation}
So we can get
$$L_{N}^{u}\geq \big| L_{n}^{u}-2L_{2n}^{u}\big|-C_{0}SnN^{-1}\geq L_{n}^{u}-2\big|L_{n}^{u}-L_{2n}^{u}\big|-C_{0}SnN^{-1}\geq \gamma S-2\big|L_{n}^{u}-L_{2n}^{u}\big|- C_{0}SnN^{-1}.$$
To obtain the second inequality (\ref{eq17}), observe that by virtue of all arguments so far apply equally well to $M_{2N}^{u}$ instead of $M_{N}^{u}$, that is
$$\big|L_{2N}^{u}+L_{n}^{u}-2L_{2n}^{u} \big|\leq C_{0}SnN^{-1},$$
so $$L_{N}^{u}-L_{2N}^{u}\leq \big|L_{N}^{u}+L_{n}^{u}-2L_{2n}^{u} \big|+\big|L_{2N}^{u}+L_{n}^{u}-2L_{2n}^{u} \big|\leq C_{0}SnN^{-1}.$$
Denote
$$u_{N}^{u}(x,y)=\frac{1}{N}\log\lVert M_{N}^{u}(x,y)\rVert,\ u_{n}^{u}(x,y)=\frac{1}{n}\log\lVert M_{n}^{u}(x,y)\rVert,\ u_{2n}^{u}(x,y)=\frac{1}{2n}\log\lVert M_{2n}^{u}(x,y)\rVert $$
In view of (\ref{eq55}) (\ref{eq10}), both $u_{n}^{u}$ and $u_{2n}^{u}$ extends to separately subharmonic functions in both variables satisfying
$$u_{n}^{u}(z_{1},z_{2})\leq S(\lambda,E),\ u_{2n}^{u}(z_{1},z_{2})\leq S(\lambda,E).$$
Applying Lemma \ref{lemma} to $\frac{u_{n}^{u}}{S}$ and $\frac{u_{2n}^{u}}{2S}$
$$  \textrm{mes}\Big[(x,y)\in \mathbb{T}^{2}:\Big|\frac{1}{N}\sum_{k=1}^{N} \frac{u_{n}^{u}}{S}\circ T^{k}_{\omega}(x,y)-\frac{1}{S}\langle u_{n}^{u}\rangle\Big|>N^{-\frac{1}{12}+2\delta}\Big]\leq C \exp(-cN^{\delta}),$$
$$\textrm{mes}\Big[(x,y)\in \mathbb{T}^{2}:\Big|\frac{1}{N}\sum_{k=1}^{N} \frac{u_{2n}^{u}}{S}\circ T^{k}_{\omega}(x,y)-\frac{1}{S}\langle u_{2n}^{u}\rangle\Big|>N^{-\frac{1}{12}+2\delta}\Big]\leq C \exp(-cN^{\delta}).$$
So there is a set $\mathcal{B}_{2}\subset\mathbb{T}^{2}$ with measure
\begin{equation}\textrm{mes(}\mathcal{B}_{2})\leq C\exp(-N^{\delta}),\end{equation}
such that for any $(x,y)\in \mathcal{G}:=\mathbb{T}^{2}\setminus(\mathcal{B}_{1}\cup\mathcal{B}_{2})$
\begin{align}\label{delta}
&\Big|  \frac{1}{N}\log\lVert M_{N}^{u}(x,y)\rVert+L_{n}^{u}-2L_{2n}^{u}  \Big|\notag\\
&\leq\Big|\frac 1N\log\lVert M_{N}^{u}(x,y)\rVert +\frac{1}{N}\sum_{k=1}^{N} \frac 1n\log\lVert M_{n}^{u}\circ T^{k}(x,y)\rVert -\frac{2}{N}\sum_{k=1}^{N}\frac{1}{2n}\log\lVert M_{2n}^{u}\circ T^{k}(x,y)\rVert \Big|\notag\\
&+\Big|\frac{1}{N}\sum_{k=1}^{N}\frac{1}{n}\log \lVert M_{n}^{u}\circ T^{k}(x,y)\rVert- L_{n}^{u} \Big|+\Big| \frac{2}{N}\sum_{k=1}^{N}\frac{1}{2n}\log\lVert M_{2n}^{u}\circ T^{k}(x,y)\rVert-2L_{2n}^{u}  \Big|\notag\\
&\leq CSnN^{-1}+C_{\delta}SN^{-\frac{1}{12}+2\delta}.
\end{align}

For small $\delta$ the second term of (\ref{delta}) is the larger one since $N\geq n^{2}$. Fix such an integer $N$. Consider the following decomposition of $u^{u}:=u_{N}^{u}$ as a function of $\mathbb{T}^{2}$:
$$u^{u}=u^{u}\chi_{\mathcal{G}}+L_{N}^{u}\chi_{\mathcal{G}^{c}}+u^{u}\chi_{\mathcal{G}^{c}}-L_{N}^{u}\chi_{\mathcal{G}^{c}}:=u_{0}^{u}+u_{1}^{u}.$$
Here $u_{0}^{u}$ is the sum of the first two terms (and $\mathcal{G}^{c}:=\mathbb{T}^{2}\setminus\mathcal{G}$)
\begin{align}
\lVert u_{0}^{u}-\langle u^{u}\rangle\rVert_{\infty}=\lVert u_{0}^{u}-L_{N}^{u}\rVert_{\infty}&=\lVert u^{u}-L_{N}^{u}\rVert_{L^{\infty}(\mathcal{G})}\notag\\
&\leq\lVert u_{N}^{u}+L_{n}^{u}-L_{2n}^{u} \rVert_{L^{\infty}(\mathcal{G})}+|L_{N}^{u}+L_{n}^{u}-L_{2n}^{u}|\notag\\
&\leq C_{\delta}SN^{-\frac{1}{12}+2\delta}.
\end{align}
On the other hand, we have
\begin{equation}\lVert u_{1}^{u} \rVert_{L^{1}(\mathbb{T}^{2})}=\int_{\mathcal{G}^{c}}|u^{u}-L_{N}^{u}|dxdy\leq2S mes(\mathcal{G}^{c})\leq CS[N^{-9}+\exp(-N^{\delta})]\leq C_{\delta}SN^{-9}.\end{equation}
Applying lemma 3.2. to the function $\frac{u^{u}}{S}$, with $\varepsilon_{0}=C_{\delta}N^{-\frac{1}{12}+2\delta}$ and $\varepsilon_{1}=C_{\delta}N^{-9}$
\begin{equation}\textrm{mes}\Big[(x,y)\in\mathbb{T}^{2}: \Big|\frac{u_{N}^{u}}{S}-\frac{1}{S}L_{N}^{u}\Big|>B^{\delta}  \log\frac{N}{\varepsilon_{1}}\Big]\leq CN^{2}\varepsilon_{1}^{-1}\exp(-cB^{-\frac 12+\delta}),\end{equation}
where $B=\varepsilon_{0}\log\frac{N}{\varepsilon_{1}}+N^{\frac{3}{2}}\varepsilon_{1}^{\frac 14}=C_{\delta}N^{-\frac {1} {12}+2\delta}\log(C_{\delta}^{-1}N^{10})+N^{\frac{3}{2}}C_{\delta}^{\frac 14}N^{-\frac{9}{4}}=N^{-\tau^{\prime}},$
$$B^{\delta}\log\frac{N}{\varepsilon_{1}}=N^{-\tau^{\prime}\delta}\log( C_{\delta}^{-1}N^{10})=N^{-\tau},$$
$$CN^{2}\varepsilon_{1}^{-1}\exp(-cB^{-\frac 12+\delta})=CN^{2}C_{\delta}^{-1}N^{9}\exp(-cN^{\tau^{\prime}(\frac 12-\delta)}) \leq C\exp(-N^{\sigma}),$$
here $\sigma\ll\delta$, so when $\delta<\frac{1}{24}$, we can get the large deviation theorem
\begin{equation}\textrm{mes}\Big[(x,y)\in\mathbb{T}^{2}: \Big|\frac{1}{N}\log\lVert M_{N}^{u}(x,y;\lambda,E)\rVert-L_{N}^{u}(\lambda,E)\Big|>S(\lambda,E)N^{-\tau}]\leq C\exp(-N^{\sigma}).\end{equation}
\begin{lemma}{\label{leite}} There are large constants $\lambda_{0}$, and $B$ depending on $v$ and $a$, such that for any positive integer $n$,
 \begin{equation}\label{eq20}
 \sup_{E}\textrm{mes}\Big[(x,y)\in\mathbb{T}^{2}:\Big|\frac{1}{n}\log\lVert M_{n}^{u}(x,y;\lambda,E)\rVert-L_{n}^{u}(\lambda,E)\Big|> \frac{1}{20}S(\lambda,E)\Big]<n^{-50},
 \end{equation}
provided $\lambda\geq \lambda_{0}\vee n^{B}$. Further for those $\lambda$ and for all $E$\\
$$L_{n}^{u}(\lambda,E)\geq\frac{1}{2}S(\lambda,E)\;\;and \;\;L_{n}^{u}(\lambda,E)-L_{2n}^{u}(\lambda,E)\leq\frac{1}{80}S(\lambda,E).$$
\end{lemma}
\Proof We have got the fundamental matrix of the Jacobi operator
$$ M_{n}(x,y;\lambda,E)=\Bigg(\begin{array}{ccc}   \frac{f_{n}(x,y;\lambda,E)}{\prod_{j=2}^{n+1}a_{j}} & -\frac{a_{1}}
{a_{2}}\frac{f_{n-1}(T(x,y);\lambda,E)}{\prod_{j=3}^{n+1}a_{j}}\\
\frac{ f_{n-1}(x,y;\lambda,E)}{\prod_{j=2}^{n}a_{j}}&-\frac{a_{1}}{a_{2}}\frac{f_{n-2}(T(x,y);\lambda,E)}
{\prod_{j=3}^{n}a_{j}}\end{array} \Bigg)$$
here
\[f_{n}(x,y;\lambda,E) =\det \left( {\begin{array}{*{20}{c}}
\lambda v_{1}-E&-a_{2}&0& 0&\cdots&\cdots &0\\
-a_{2}&\lambda v_{2}-E &-a_{3}&0&\cdots&\cdots&0\\
0&-a_{3}& \lambda v_{3}-E&-a_{4}&\cdots&\cdots&0\\
\vdots & \vdots &\vdots &\vdots &{}& \ddots &-a_{n}\\
0&0&{0}&{\cdots}&0&-a_{n}&\lambda v_{n}-E
\end{array}} \right)\]
The matrix on the right-hand side can be written in the form $D_{n}+B_{n}$
$$D_{n}(x,y;\lambda,E)=\textrm{diag}(\lambda v_{1}-E,...,\lambda v_{n}-E).$$
Because $1\leq|a_{n}|\leq2$, clearly, $ \lVert B_{n}  \rVert \leq 4$ and
\begin{equation}\frac{1}{n}\log|\det D_{n}(x,y;\lambda,E)|=\log\lambda+\frac{1}{n}\sum_{j=1}^{n}\log\big|v(T^{j}_{\omega}(x,y))-E/\lambda\big|.\end{equation}

Because of the classical Lojasiewicz result, for nonconstant real-analytic functions $v$, there exist constant $b>0$ and C depending on $v$ such that
\begin{equation}\label{mes}
\textrm{mes}[(x,y)\in\mathbb{T}^{2}:|v(x,y)-h|<t]\leq Ct^{b},
\end{equation}
for any $-2\lVert v\rVert_{\infty}\leq |h|\leq 2\lVert v\rVert_{\infty}$ and $t>0$.
We divide $E$ into two cases.\\

\textbf{case 1.} When $E\leq 2\lambda\lVert v\rVert_{\infty}$,
\begin{equation}\textrm{mes}\Big[(x,y)\in \mathbb{T}^{2}:\frac{1}{n}\sum_{j=1}^{n}\log|v\circ T^{j}_{\omega}(x,y)-E/\lambda|\leq -\rho\Big]\leq nCe^{-b\rho},\end{equation}
one also have the upper bound
\begin{equation}\sup_{(x,y)\in\mathbb{T}^{2}}\frac{1}{n}\sum_{j=1}^{n}\log|v\circ T^{j}_{\omega}(x,y)-E/\lambda|\leq \log(3\lVert v\rVert_{\infty}).\end{equation}
Since $$\lVert D_{n}^{-1}(x,y;\lambda,E)\rVert\leq \lambda^{-1}\max_{1\leq j\leq n}|v\circ T^{j}_{\omega}(x,y)-E/\lambda|^{-1},$$
(\ref{mes}) implies that \begin{align}
&\textrm{mes}[(x,y)\in\mathbb{T}^{2}:\lVert D_{n}^{-1}(x,y;\lambda,E)\rVert > \frac{1}{8}]\notag\\
&\leq n \;mes[(x,y)\in\mathbb{T}^{2}:|v(x,y)-E/\lambda|<8\lambda^{-1}]\notag\\
&\leq Cn\lambda^{-b}.
\end{align}
Hene \begin{equation}\textrm{mes}[(x,y)\in\mathbb{T}^{2}:\lVert D_{n}^{-1}(x,y;\lambda,E)B_{n}\rVert>\frac{1}{2}]\leq Cn\lambda^{-b}.\end{equation}
So we have
\begin{align}\label{eq18}
\Big|\frac{1}{n}\log\big|\frac{f_{n}(x,y;\lambda,E)}{\prod_{j=2}^{n+1}a_{j}}\big|-\log \lambda\Big|
&\leq\big|\frac{1}{n}\log|f_{n}(x,y;\lambda,E)|-\log \lambda\big|+\big|\frac 1n\log \prod_{j=2}^{n+1}a_{j}\big|\notag\\
&\leq\big|\frac{1}{n}\log|f_{n}(x,y;\lambda,E)|-\log \lambda\big|+\log2\notag\\
&\leq\Big|\frac{1}{n}\sum_{j=1}^{n}\log|v\circ T^{j}(x,y)-E/\lambda|\Big|+\Big|\frac{1}{n}\log|\det(I+D_{n}(x,y;\lambda,E)^{-1}B_{n})|\Big|+\log2\notag\\
&\leq \rho +\log(3\lVert v\rVert_{\infty})+2\log 2,
\end{align}
up to a set of measure not exceeding
\begin{equation}\label{eq19}
Cne^{-b\rho}+Cn\lambda^{-b}.
\end{equation}

Now let $\rho=\frac{1}{400}\log \lambda$ and assume that $\lambda\geq(12\lVert v\rVert_{\infty})^{400}$. the right-hand side of (\ref{eq18}) is no larger than $\frac{1}{200}\log \lambda$. The first part of (\ref{eq19}) is the larger one under these assumptions, choosing large constant $B$ depending only on $v$ such that when $\lambda\geq n^{B}$, we have
$$Cne^{-b\rho}+Cn\lambda^{-b}\leq Cn\lambda^{-\frac{b}{400}}\leq n^{-100}$$
\begin{equation}\sup_{|E|\leq2\lVert v\rVert_{\infty}}\textrm{mes}\Big[(x,y)\in\mathbb{T}^{2}: \Big|\frac{1}{n}\log\big|\frac{f_{n}(x,y;\lambda,E)}{\prod_{j=2}^{n+1}a_{j}}\big|-\log \lambda\Big|\geq\frac{1}{200}\log\lambda\Big]\leq n^{-100}.\end{equation}
Similarly, we have
\begin{equation}\sup_{|E|\leq2\lVert v\rVert_{\infty}}\textrm{mes}\Big[(x,y)\in\mathbb{T}^{2}: \Big|\frac{1}{n}\log\big|\frac{f_{n-1}(x,y;\lambda,E)}{\prod_{j=2}^{n}a_{j}}\big|-\log \lambda\Big|\geq\frac{1}{200}\log\lambda       \Big]\leq (n-1)^{-100}.\end{equation}
In view of the relationship between $M_{n}(x,y;\lambda,E)$ and $f_{n}(x,y;\lambda,E)$, one therefore obtains
$$\sup_{|E|\leq2\lVert v\rVert_{\infty}}\textrm{mes}\Big[(x,y)\in\mathbb{T}^{2}: \Big|\frac{1}{n}\log\lVert M_{n}(x,y;\lambda,E)\rVert-\log \lambda\Big|\geq\frac{1}{199}\log\lambda       \Big]\leq 4(n-2)^{-100}.$$
Because of the relationship between $\frac{1}{n}\log\lVert M_{n}(x,y;\lambda,E)\rVert$ and $\frac{1}{n}\log\lVert M_{n}^{u}(x,y;\lambda,E)\rVert$, we have $$\sup_{|E|\leq2\lVert v\rVert_{\infty}}\textrm{mes}\Big[(x,y)\in\mathbb{T}^{2}: \Big|\frac{1}{n}\log\lVert M_{n}^{u}(x,y;\lambda,E)\rVert-\log \lambda\Big|\geq\frac{1}{198}\log\lambda    \Big]\leq 4(n-2)^{-100},$$
In particular,
\begin{equation}|L_{n}^{u}(\lambda,E)-\log\lambda|\leq \frac{1}{198}\log\lambda+4S(\lambda,E)(n-2)^{-100}\leq \frac{1}{190}S(\lambda,E),\end{equation}
provided $n> 2$. Since
$$\log \lambda\geq\frac{99}{100}\sup_{|E|\leq2\lambda\lVert v\rVert_{\infty}}S(\lambda,E)$$
for large $\lambda_{0}$.
So, we have
\begin{equation}
L_{n}^{u}(\lambda,E)\geq\log\lambda-\frac{1}{190}S(\lambda,E)\geq\frac{99}{100}\sup_{|E|\leq2\lambda\lVert v\rVert_{\infty}}S(\lambda,E)-\frac{1}{190}S(\lambda,E)\geq\frac{1}{2}S(\lambda,E),\end{equation} \begin{equation}L_{n}^{u}(\lambda,E)-L_{2n}^{u}(\lambda,E)\leq|L_{n}^{u}(\lambda,E)-\log\lambda|+|L_{2n}^{u}(\lambda,E)-\log\lambda|
\leq\frac{1}{80}S(\lambda,E)\end{equation}
and
\begin{align}
&\sup_{|E|\leq2\lVert v\rVert_{\infty}}\textrm{mes}\Big[(x,y)\in\mathbb{T}^{2}: \Big|\frac{1}{n}\log\lVert M_{n}^{u}(x,y;\lambda,E)\rVert-L_{n}^{u}(\lambda,E)\Big|\geq\frac{1}{20}S(\lambda,E)\Big]\notag\\
&\leq \sup_{|E|\leq2\lVert v\rVert_{\infty}}\textrm{mes}\Big[(x,y)\in\mathbb{T}^{2}: \Big|\frac{1}{n}\log\lVert M_{n}^{u}(x,y;\lambda,E)\rVert-\log \lambda\Big|+|L_{n}^{u}(\lambda,E)-\log\lambda|\geq\frac{1}{20}S(\lambda,E)\Big]\notag\\
&\leq\sup_{|E|\leq2\lVert v\rVert_{\infty}}\textrm{mes}\Big[(x,y)\in\mathbb{T}^{2}: \Big|\frac{1}{n}\log\lVert M_{n}^{u}(x,y;\lambda,E)\rVert-\log \lambda\Big|\geq (\frac{1}{20}-\frac{1}{190})S(\lambda,E)\Big]\notag\\
&\leq\sup_{|E|\leq2\lVert v\rVert_{\infty}}\textrm{mes}\Big[(x,y)\in\mathbb{T}^{2}: \Big|\frac{1}{n}\log\lVert M_{n}^{u}(x,y;\lambda,E)\rVert-\log \lambda\Big|\geq\frac{1}{198}\log\lambda  \Big]\notag\\
&\leq 4(n-2)^{-100}\leq n^{-50}.
\end{align}

\textbf{case 2.} When $|E|>2\lambda \lVert v\rVert_{\infty}$ and $\lambda_{0}$ is sufficiently large, then the set in (\ref{eq20}) is empty. In fact, for such $E$,
$$\Big|\frac{1}{n}\log|\det D_{n}(x,y;\lambda,E)|-\log|E|\Big|=\Big|\frac 1n\sum_{j=1}^{n}\log\big|\frac{\lambda v(T_{\omega}^{j})-E}{E}\big| \Big| \leq 2,$$
and thus
$$\Big| \frac{1}{n}\log|f_{n}(x,y;\lambda,E)|-\log|E|\Big|\leq 4.$$
So$$\Big| \frac{1}{n}\log\big|\frac{f_{n}(x,y;\lambda,E)}{\prod_{j=2}^{n+1}a_{j}}\big|-\log|E|\Big|\leq 4+\log2$$
and $$\Big| \frac{1}{n}\log\big|\frac{f_{n-1}(x,y;\lambda,E)}{\prod_{j=2}^{n}a_{j}}\big|-\log|E|\Big|\leq 4+\log2,$$
which implies that for large $\lambda$,
$$\Big| \frac{1}{n}\log\lVert M_{n}(x,y;\lambda,E)\rVert-\log|E|\Big|\leq8+2\log2\leq\frac{1}{400}S(\lambda,E),$$
$$\Big| \frac{1}{n}\log\lVert M_{n}^{u}(x,y;\lambda,E)\rVert-\log|E|\Big|\leq\frac{1}{400}S(\lambda,E)+ \frac12\log2\leq\frac{1}{200}S(\lambda,E).$$
So
$$\Big| L_{n}^{u}(\lambda,E)-\log|E|\Big|\leq\frac{1}{200}S(\lambda,E),$$
\begin{equation}
\big|\frac{1}{n}\log\lVert M_{n}^{u}(x,y;\lambda,E)\rVert-L_{n}^{u}(\lambda,E)\big|\leq \frac{1}{200}S(\lambda,E)+\frac{1}{200}S(\lambda,E)\leq\frac{1}{20}S(\lambda,E)\end{equation}
and
\begin{align}
L_{n}^{u}(\lambda,E)\geq\log|E|-\frac{1}{200}S(\lambda,E)&\geq\log(2\lambda\lVert v\rVert_{\infty})-\frac{1}{200}S(\lambda,E)\notag\\
&\geq\frac{99}{100}\sup_{|E|\geq2\lambda\lVert v\rVert_{\infty}}S(\lambda,E)-\frac{1}{200}S(\lambda,E)\notag\\
&\geq\frac{1}{2}S(\lambda,E),
\end{align}
\begin{equation}L_{n}^{u}(\lambda,E)-L_{2n}^{u}(\lambda,E)\leq\Big| L_{n}^{u}(\lambda,E)-\log|E|\Big|+\Big| L_{2n}^{u}(\lambda,E)-\log|E|\Big|\leq\frac{1}{80}S(\lambda,E),\end{equation}
and the lemma follows. \qquad \qedbox
\section{The large deviation theorem and continuity of the Lyapunov exponent}
\noindent \textbf{Proof of Theorem \ref{theorem1}.}

Fix $\sigma<\frac{1}{24}$ throughout the proof and let $\tau=\tau(\sigma)>0$ be as in (\ref{eq21}). Moreover, let $\lambda\geq \lambda_{0}\vee n^{B}:=\lambda_{1}$ be as in Lemma \ref{leite}, in this proof we shall require $n_{0}$ to be sufficiently large at various places, but of course $n_{0}$ will be assumed fixed. In view of Lemma \ref{leite} the hypothesis of Lemma \ref{lemma1} are satisfied with $\gamma=\gamma_{0}=\frac 12$,
\begin{equation}\label{eq22}
n_{0}^{2}<N<n_{0}^{5},
\end{equation}
provided
\begin{equation}
9n_{0}\geq 20\log(2n_{0}^{5}).
\end{equation}
It is clear that holds if $n_{0}$ is large. Applying Lemma \ref{lemma1}. one obtains (suppressing $\lambda$, E for simplicity)
\begin{equation}
L_{N}^{u}\geq(\frac 12- \frac{1}{40})S-C_{0}SN^{-1}n_{0}\geq \gamma_{1}S,
\end{equation}
\begin{equation}
L_{N}^{u}-L_{2N}^{u}\leq C_{0}SN^{-1}n_{0}\leq\frac{\gamma_{1}}{40}S,
\end{equation}
with $\gamma_{1}=\frac 13$. Moreover, with some constant $C_{1}\geq 1$ depending on $\varepsilon$,
\begin{equation}\label{eq23}
\textrm{mes}\Big[(x,y)\in\mathbb{T}^{2}:\big|\frac{1}{N}\log\lVert M_{N}^{u}(x,y;\lambda,E)\rVert-L_{N}^{u}(\lambda,E)\big|> S(\lambda,E)N^{-\tau}\Big]<C_{1}\exp(-N^{\sigma}),
\end{equation}
for all $N$ in the range given by (\ref{eq22}).

So
\begin{equation}\label{eq29}
\textrm{mes}\Big[(x,y)\in\mathbb{T}^{2}:\big|\frac{1}{N}\log\lVert M_{N}(x,y;\lambda,E)\rVert-L_{N}(\lambda,E)\big|> S(\lambda,E)N^{-\tau}\Big]<C_{1}\exp(-N^{\sigma}),
\end{equation}
in particular, (\ref{eq23}) implies that
$$\textrm{mes}\Big[(x,y)\in\mathbb{T}^{2}:\big|\frac{1}{N}\log\lVert M_{N}^{u}(x,y;\lambda,E)\rVert-L_{N}^{u}(\lambda,E)\big|> S(\lambda,E)\frac{\gamma_{1}}{10}\Big]<C_{1}\exp(-N^{\sigma})\leq \bar{N}^{-10},$$
provided $n_{0}$ is large and
$$N^{2}\leq\bar{N}\leq C_{1}^{-\frac{1}{10}}\exp(\frac{1}{10} N^{\sigma}).$$
The first inequality was added to satisfy (\ref{eq24}). In view of (\ref{eq22}), one thus has the range
\begin{equation}\label{eq25}
n_{0}^{4}\leq\bar{N}\leq \exp(\frac{1}{10}n_{0}^{5\sigma})
\end{equation}
of admissible $\bar{N}$. Moreover,
\begin{equation}\label{eq26}L_{\bar{N}}^{u}\geq\gamma_{1}S-2C_{0}SN^{-1}n_{0}-C_{0}S\bar{N}^{-1}N\end{equation}
and \begin{equation}\label{eq27}L_{\bar{N}}^{u}-L_{2\bar{N}}^{u}\leq C_{0}S\bar{N}^{-1}N.\end{equation}

At the next stage of this procedure, observe that the left end-point of the range of admissible indices starts at $n_{0}^{8}$, which is less than the right end-point of the range (\ref{eq25}) (for $n_{0}$ large). Therefore, from this point on the ranges will overlap and cover all large integers. To ensure that the process does not terminate, simply note the rapid convergence of the series given by (\ref{eq26}) (\ref{eq27}).
so \begin{equation}\label{eq28}L(\lambda,E)=L^{u}(\lambda,E)=\inf\limits_{n}L_{n}^{u}(\lambda,E)\geq \frac{1}{4}\log \lambda.\end{equation}
Thus the theorem follows because of (\ref{eq29}) (\ref{eq28}).\qquad \qedbox\\

\noindent \textbf{Proof of Theorem \ref{theorem2}.}

Let $\lambda\geq \lambda_{0}\vee n^{B}:=\lambda_{1}$ be as in Lemma \ref{leite}, fix any positive $\sigma<\frac{1}{24}$. And we set $n=\lfloor C_{0}(\log N)^{\frac{1}{\sigma}}\rfloor$ with some large constant $C_{0}$. For large integer $n$,
$$L_{n}^{u}(E)\geq\frac{1}{2}S(\lambda,E)\;\;and \;\;L_{n}^{u}(E)-L_{2n}^{u}(E)\leq\frac{1}{80}S(\lambda,E).$$
 Let $m=\frac{N}{n}$, then by (\ref{eq31}), we have for $0\leq j\leq m-1$
$$\Big|\frac{1}{n}\log\lVert M_{n}^{u}(T^{jn}_{\omega}(x,y);\lambda,E)\rVert-L_{n}^{u}(E)\Big|\leq S(\lambda,E)n^{-\tau}$$
$$\Big|\frac{1}{2n}\log\lVert M_{2n}^{u}(T^{jn}_{\omega}(x,y);\lambda,E)\rVert-L_{2n}^{u}(E)\Big| \leq S(\lambda,E)n^{-\tau}$$
for $(x,y)\in \mathbb{G}$, with
$$\textrm{mes}(\mathbb{T}^{2}\backslash\mathbb{G})\leq 2m\exp(-cn^{\sigma})\leq CN\exp(-n^{\sigma})\leq \frac{1}{N^{2}}$$
Thus when $(x,y)\in\mathbb{G}$,
\begin{equation}\lVert M_{n}^{u}(T^{jn}_{\omega}(x,y);\lambda,E)\rVert>\exp(n(L_{n}^{u}(E)-Sn^{-\tau}))>\exp(\frac 12 nS-Snn^{-\tau})>\exp(\frac{1}{4}nS):=\mu,
\end{equation}
and
\begin{align}&\log\lVert M_{n}^{u}(T^{jn}_{\omega}(x,y);\lambda,E)\rVert+\log\lVert M_{n}^{u}(T^{(j+1)n}_{\omega}(x,y);\lambda,E)\rVert-\log\lVert M_{n}^{u}(T^{(j+1)n}_{\omega}(x,y);\lambda,E) M_{n}^{u}(T^{jn}_{\omega}(x,y);\lambda,E)\rVert\notag\\
&\leq 2n(L_{n}^{u}(E)+Sn^{-\tau})-2n(L_{2n}^{u}(E)-Sn^{-\tau})\notag\\
&\leq 4Sn^{1-\tau}+2n(L_{n}^{u}(E)-L_{2n}^{u}(E))\notag\\
&\leq 4Sn^{1-\tau}+\frac{1}{40}nS\leq \frac 12\log\mu.
\end{align}
By applying the avalanche principle \ref{Avalanche}, and integrating over $\mathbb{G}$ one obtains
\begin{align}\label{eq60}
&\Big|\int_{\mathbb{G}}u_{N}^{u}(x,y;\lambda,E)dxdy+\frac{1}{m}\int_{\mathbb{G}}\sum_{j=2}^{m-1}u_{n}^{u}(T_{\omega}^{(j-1)n}(x,y);\lambda,E)dxdy
-\frac{2}{m}\int_{\mathbb{G}}\sum_{j=1}^{m-1}u_{2n}^{u}(T_{\omega}^{(j-1)n}(x,y);\lambda,E)dxdy  \Big|\notag\\
&<C\frac{m}{N\mu}\leq \frac{Cn}{N}
\end{align}
We want to replace here the integration over $\mathbb{G}$ by integration over $\mathbb{T}^{2}$, recall that due to (4) in Remark 1.1
$$\int_{\mathbb{T}^{2}}(u_{n}^{u}(x,y;\lambda,E))^{2}dxdy\leq \tilde{C}(\lambda,v,a)$$
for any $n$ and any $E$. Here, by Cauchy-Schwartz inequality
$$|\int_{\mathbb{B}}u_{n}^{u}(x,y;\lambda,E)dxdy|\leq\tilde{C}(\lambda,v,a)^{\frac 12}(\textrm{mes} \mathbb{B})^{\frac 12}$$
for any $n$, any $E$ and any $\mathbb{B}\subseteq \mathbb{T}^{2}$. Thus
\begin{align}\label{eq61}
&\Big|\int_{\mathbb{T}^{2}\backslash\mathbb{G}}u_{N}^{u}(x,y;\lambda,E)dxdy+\frac{1}{m}\int_{\mathbb{T}^{2}\backslash\mathbb{G}}\sum_{j=2}^{m-1}u_{n}^{u}(T_{\omega}^{(j-1)n}(x,y);\lambda,E)dxdy
-\frac{2}{m}\int_{\mathbb{T}^{2}\backslash\mathbb{G}}\sum_{j=1}^{m-1}u_{2n}^{u}(T_{\omega}^{(j-1)n}(x,y);\lambda,E)dxdy  \Big|\notag\\
&\leq 4\tilde{C}(\lambda,v,a)^{\frac 12}(\textrm{mes}(\mathbb{T}^{2}\backslash\mathbb{G}))^{\frac12}\leq \frac{Cn}{N}
\end{align}
Combining \ref{eq60} with \ref{eq61}, we have
$$\Big| L_{N}^{u}(E)+\frac{m-2}{m}L_{n}^{u}(E)-\frac{2(m-1)}{m}L_{2n}^{u}(E)\Big|\leq \frac{Cn}{N}$$
Thus
\begin{align}\label{eq33}
|L_{N}^{u}(E)-2L_{2n}^{u}(E)+L_{n}^{u}(E)|&\leq \frac{Cn}{N}+\frac{2}{m}|L_{n}^{u}(E)-L_{2n}^{u}(E)|\notag\\
&\leq\frac{Cn}{N}+\frac{4n}{N}C^{\prime\prime}(\lambda,v,a)\leq \frac{Cn}{N}.
\end{align}
Similarly
\begin{equation}\label{eq34}
|L_{2N}^{u}(E)-2L_{2n}^{u}(E)+L_{n}^{u}(E)|\leq \frac{Cn}{N}.
\end{equation}
 Taking the difference of the two inequalities (\ref{eq33}) (\ref{eq34}), we have
\begin{equation}
|L_{N}^{u}(E)-L_{2N}^{u}(E)|\leq \frac{Cn}{N},
\end{equation}
which after summing over dyadic $N$ gives
\begin{equation}\label{eq35}
|L_{N}^{u}(E)-L^{u}(E)|\leq \frac{Cn}{N}.
\end{equation}
So by using the (\ref{eq33}) (\ref{eq35}), we can get
\begin{equation}\label{eq36}
|L^{u}(E)-2L_{2n}^{u}(E)+L_{n}^{u}(E)|\leq \frac{Cn}{N}.
\end{equation}

Fix $\lambda$, for any $n$
$$\lVert A_{n}^{\prime}(x,y;\lambda,E)\rVert\leq (E_{0}+\lambda\lVert v\rVert_{\infty}+C_{a}):=C_{v,a}.$$
Assume for instance that $\lVert M_{n}^{a}(E)\rVert>\lVert M_{n}^{a}(E^{\prime})\rVert$
\begin{align*}
 \log \lVert M_{n}^{a}(E)\rVert-\log\lVert M_{n}^{a}(E^{\prime})\rVert=\log\frac{\lVert M_{n}^{a}(E)\rVert}{\lVert M_{n}^{a}(E^{\prime})\rVert}&=\log(1+\frac{\lVert M_{n}^{a}(E)\rVert-\lVert M_{n}^{a}(E^{\prime})\rVert}{\lVert M_{n}^{a}(E^{\prime})\rVert})\\
&\leq \frac{\lVert M_{n}^{a}(E)\rVert-\lVert M_{n}^{a}(E^{\prime})\rVert}{\lVert M_{n}^{a}(E^{\prime})\rVert}.\\
\end{align*}
Obviously, for any $n$, $\lVert\partial_{E}M_{n}^{a}(E)\rVert\leq n(C_{v,a})^{n-1},$
so
\begin{equation}
|L_{n}^{u}(E)-L_{n}^{u}(E^{\prime})|=|L_{n}^{a}(E)-L_{n}^{a}(E^{\prime})|\leq C_{v,a}^{n}|E-E^{\prime}|.
\end{equation}
In view of this fact, (\ref{eq36}) implies that for any $E$, $E^{\prime}$,
\begin{align}
\Big|L(E)-L(E^{\prime})\Big|=\Big|L^{u}(E)-L^{u}(E^{\prime})\Big|&\leq \frac{Cn}{N}+C^{n}|E-E^{\prime}|\notag\\
&\leq e^{-cn^{\sigma}}+C^{n}|E-E^{\prime}|\notag\\
&\leq C\exp(-c(\log|E-E^{\prime}|^{-1})^{\sigma}),
\end{align}
and the theorem follows from $n\sim \log\frac{1}{|E-E^{\prime}|}$. \qquad \qedbox

\vskip1cm
\noindent{$\mathbf{Acknowledgments}$}

The authors were supported by the NSFC (grant no. 11571327)
\vskip1cm
\noindent{$\mathbf{References}$}

\end{document}